# Directed Transmission Method,

# A Fully Asynchronous and Distributed Algorithm

# to Solve Sparse Linear System

Fei Wei

Huazhong Yang

Department of Electronic Engineering, Tsinghua University, Beijing, China

Abstract

In this paper, we propose a new distributed algorithm, called Directed Transmission Method (DTM). DTM is a fully asynchronous and continuous-time iterative algorithm to solve SPD sparse linear system. As an architecture-aware algorithm, DTM could be freely running on all kinds of heterogeneous parallel computer. We proved that DTM is convergent by making use of the final-value theorem of Laplacian Transformation. Numerical experiments show that DTM is stable and efficient.

1. Introduction

Solving the large sparse linear system, $\mathbf{Ax} = \mathbf{b}$, is fundamental in the scientific computing. When the coefficient matrix $\mathbf{A}$ is symmetric-positive-definite (SPD), the linear system is called the SPD linear system, which is widely encountered in engineering applications [1]. To solve large scale SPD systems, Domain Decomposition Method (DDM) is frequently used . DDM could be classified into Schur Complement method, Additive Schwarz method (or block-Jacobi) and Multiplicative Schwarz method (or block-Gauss-Seidel) [2, 3, 4]. Most of these distributed numerical algorithms acquire synchronization.

Since synchronization is troublesome and time-consuming work, researchers have made many efforts for the asynchronous algorithms [17, 18, 19]. However, asynchronous methods were not considered mainstream by researchers in numerical analysis [18]. The main reason was that the performances of the traditional asynchronous algorithms, e.g. asynchronous block-Jacobi, are not comparable to the synchronous ones.

Directed transmission method (DTM) is a new asynchronous iterative numerical algorithm to solve sparse linear SPD systems [20]. It does not need any synchronization among processors when doing the distributed computing.

As a continuous-time iterative algorithm, the iterative formula of DTM is different from the traditional iterative algorithms, e.g. Gauss-Jacobi, which usually have the discrete-time iterative form, i.e. $\left[ x_1^k, x_2^k, \cdots x_n^k \right]^T = f\left( \left[ x_1^{k-1}, x_2^{k-1}, \cdots x_n^{k-1} \right]^T \right)$. In DTM, the continuous-time variable $t$ is used instead of the iterative index $k$, and the continuous-time iterative form of $\left[ x_1(t), x_2(t), \cdots x_n(t) \right]^T = f\left( \left[ x_1(t-\tau_1), x_2(t-\tau_2), \cdots x_n(t-\tau_n) \right]^T \right)$, $\tau_i \in \mathbb{R}^+$, $i = 1, 2, \cdots, n$, is used.

Here $\tau_i$ are positive real values, called the transmission delay of $x_i(t)$. Usually, $\tau_i$, $i = 1, 2, \cdots, n$, are different; if we set $\tau_1 = \tau_2 = \cdots = \tau_n = 1$, then DTM is degenerated into a discrete-time iterative algorithm, which is called Virtual Transmission Method (VTM). VTM is a special case of DTM, and we introduced VTM in [6]. DTM could be considered as a generalization of VTM.

DTM is an architecture-aware algorithm. In virtue of the directed transmission line, we may get a perfect one-to-one mapping from the transmission delay of the distributed algorithm to the communication delay of the distributed computer. This is called the Algorithm-Architecture Delay Mapping. This concept will be further described in this paper.

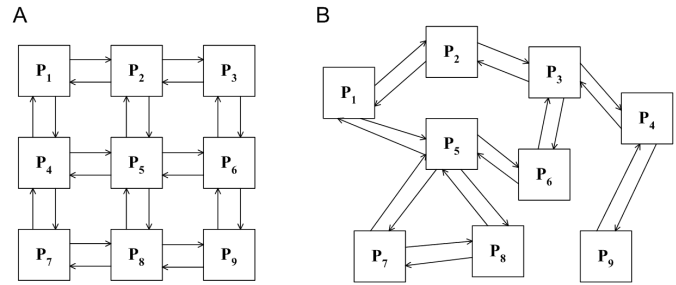

**Figure 1.** Illustration of the N2N communication model. (**A**) Regular N2N communication among processors. (**B**) Irregular N2N communication among processors.

DTM employs the Neighbor-To-Neighbor (N2N) communication model, which requires no global broadcasting but only local communication between neighboring processors, as shown in Fig. 1. The N2N model could be recognized as a kind of Peer-To-Peer (P2P) model, but it is much simpler.

DTM is inspired by the behavior of transmission lines from microwave network. The basic idea of DTM is to add Directed Transmission Lines (DTLs) into the graph of the sparse linear system to realize the asynchronous distributed computing.

This paper is organized as follows. Section 2 introduces the Directed Transmission Line (DTL). Section 3 defines the electric graph. Section 4 describes how to partition the electric graph by Electric Vertex Splitting (EVS). Section 5 details the algorithm of DTM. Section 6 presents the convergence theory for DTM. Numerical experiments are shown in Section 7. We conclude this work in Section 8.

## 2. Directed Transmission Line (DTL)

Transmission line is a magic physical element in electrical engineering [9, 10, 11]. The physical transmission line is undirected. In DTM, we bring in the directed transmission line (DTL), whose mathematical description is shown in (2.1). It should be noted that DTL does not exist in the nature, and it is an algorithmic element created by us.

(2.1) $$U_{out}(t) + Z \cdot I_{out}(t) = U_{in}(t-\tau) - Z \cdot I_{in}(t-\tau)$$

where $U_{out}(t)$ and $I_{out}(t)$ represent the potential and current of the output port, while $U_{in}(t)$ and $I_{in}(t)$ represent those of the input port. $\tau$ is the propagation delay from the input to the output. $Z$ is the characteristic impedance of DTL, which must be positive. (2.1) is called the Directed Transmission Delay Equation.

If the input and output of one DTL are the output and input of another DTL, respectively, and they have the same characteristic impedance, these two DTLs are called Directed Transmission Lines Pair (DTLP). The mathematical description of DTLP is given in (2.2).

(2.2) $$\begin{cases} U_1(t) + Z \cdot I_1(t) = U_2(t-\tau_{2\to1}) - Z \cdot I_2(t-\tau_{2\to1}) \\ U_2(t) + Z \cdot I_2(t) = U_1(t-\tau_{1\to2}) - Z \cdot I_1(t-\tau_{1\to2}) \end{cases}$$

Where $\tau_{2\to1}$ is the propagation delay of the DTL from Port 2 to Port 1, and $\tau_{1\to2}$ is the propagation delay of the DTL from Port 1 to Port 2. Fig. 2C illustrates the symbol of DTLP.

It should be noticed that $\tau_{2\to1}$ and $\tau_{1\to2}$ may be different. Further, the physical transmission line could be recognized as a special DTLP with the feature of symmetric propagation delay, i.e. $\tau_{2\to1} = \tau_{1\to2}$ [9, 10, 11]. Fig. 2 illustrates the symbols of DTL and DTLP.

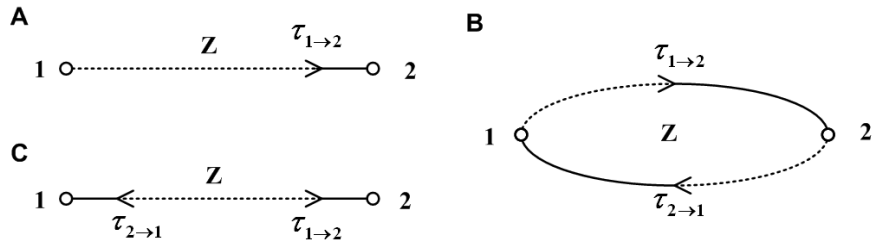

**Figure 2.** Symbols of DTL and DTLP. (**A**) The symbol of DTL from Port 1 to Port

2. **(B)** DTLP between Port 1 and Port 2. **(C)** The simplified symbol of DTLP between Port 1 and Port 2.

Why do we prefer DTL rather than the undirected one? This is because that the communication from one processor to another is directed, i.e. the communication delay from Processor A to B may be different from that from Processor B to A. In virtue of DTL, we may get a perfect one-to-one mapping from the transmission delay of DTL to the communication delay of the digital data link between processors. This is concept of Algorithm-Architecture Delay Mapping.

## 3. Electric Graph

Assume there is an *n*-dimension symmetric linear system,

(3.1) $$\mathbf{Ax} = \mathbf{b}$$

where $\mathbf{x} = (x_1, \cdots, x_n)^T$, $\mathbf{b} = (b_1, \cdots b_n)^T$, $\mathbf{A}$ is an *n*-dimension symmetric matrix.

According to the graph theory of matrix, $\mathbf{A}$ could be represented by an undirected graph $G$ [2, 3]. Assume $i \neq j$, if $a_{ij} \neq 0$, there is an edge $E_{ij}$ between $V_i$ and $V_j$ in the $G$; otherwise, $V_i$ and $V_j$ are not connected.

In this section, we define the electric graph $G_e$ of the symmetric linear system (3.1). We call $a_{ij}$ the weight of the edge $E_{ij}$, if $i \neq j$. We call $a_{ii}$ the weight of $V_i$, and call $b_i$ the source of $V_i$. $x_i$ is called the potential of $V_i$. It is easy to know that an electric graph is one-to-one mapped to a symmetric linear system. $G_e$ is defined to be SPD, iff the coefficient matrix $\mathbf{A}$ is SPD. All these concepts and terminology will be useful to describe the DTM algorithm in the following sections.

**Example 3.1**: The electric graph of (3.2) is shown in Fig. 3.

(3.2) $$\begin{pmatrix} 5 & -1 & -1 & 0 \\ -1 & 6 & -2 & -1 \\ -1 & -2 & 7 & -2 \\ 0 & -1 & -2 & 8 \end{pmatrix} \begin{pmatrix} x_1 \\ x_2 \\ x_3 \\ x_4 \end{pmatrix} = \begin{pmatrix} 1 \\ 2 \\ 3 \\ 4 \end{pmatrix}$$

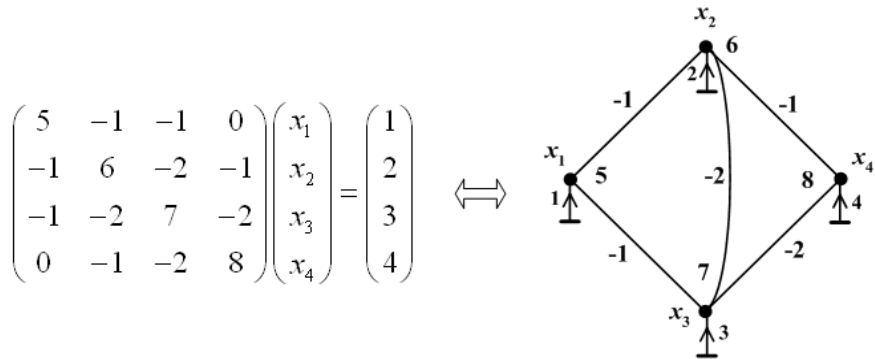

**Figure 3.** Illustration of the electric graph of the symmetric linear system (3.2).

## 4. Electric Vertex Splitting (EVS)

In this section, we partition the electric graph $G_e$ into a number of separated subgraphs by a new splitting technique, called Electric Vertex Splitting (EVS), which is also presented in [6]. EVS is called wire tearing when partitioning the circuit.

The basic idea of this partitioning technique is based on the Kirchhoff's Current Law from circuit theory [12, 16, 17, 18, 19]. This concept is illustrated in Fig. 4.

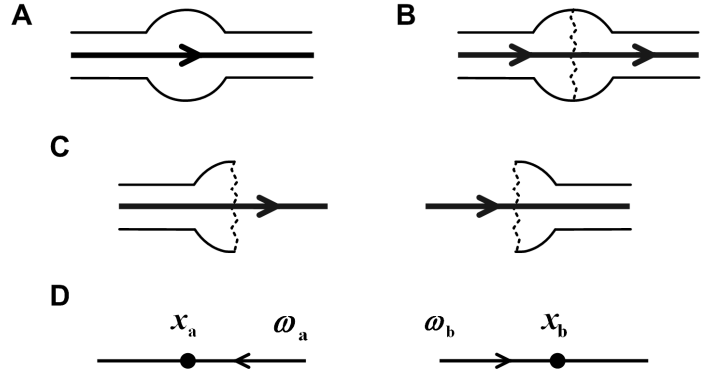

**Figure 4.** Illustration of EVS. (**A**) The original node, with current flowing through it. (**B**) Splitting this node. (**C**) The node is split into a pair of twin nodes, and the currents are disclosed. (**D**) Simplified symbol of the inflow currents.

There are four steps to do EVS.

Step-1, Set the splitting boundary $G_B$. $V \in G_e$ is called boundary vertex iff $V \in G_B$; otherwise, $V$ is called inner vertex.

Step-2, Split each boundary vertex into two vertices, which are called twin vertices.

Step-3, Split the weight and source of each boundary vertex, and split the weight of edge on the boundary, i.e. $E_{ij}$, if $E_{ij} \in G_e$ and $V_i, V_j \in G_B$.

Step-4, Add new variables, called inflow currents, $\omega$, to represent the influence coming from the adjacent subgraphs.

After these four steps, the original electric graph is split into $N$ subgraphs. If there is inflow current flowing into one vertex, then this vertex is called a port. As the result, twin vertices are also the ports of subgraphs.

**Example 4.1**: We split the electric graph $G_e$ of the linear system (3.2). $V_2$ and $V_3$ are set to be the boundary $G_B$ and we split the weights and sources of them, then we get 4 ports, $P_{2a}$, $P_{2b}$, $P_{3a}$ and $P_{3b}$, with currents $\omega_{2a}$, $\omega_{2b}$, $\omega_{3a}$ and $\omega_{3b}$ flowing into them, respectively. After that $G_e$ is split into two subgraphs. Finally we obtain two subsystems (4.1) and (4.2).

(4.1)
$$\begin{pmatrix} 5 & -1 & -1 \\ -1 & 2.5 & -0.9 \\ -1 & -0.9 & 3.3 \end{pmatrix} \begin{pmatrix} x_1 \\ x_{2a} \\ x_{3a} \end{pmatrix} = \begin{pmatrix} 1 \\ 0.8 \\ 1.6 \end{pmatrix} + \begin{pmatrix} 0 \\ \omega_{2a} \\ \omega_{3a} \end{pmatrix}$$

(4.2)
$$\begin{pmatrix} 3.5 & -1.1 & -1 \\ -1.1 & 3.7 & -2 \\ -1 & -2 & 8 \end{pmatrix} \begin{pmatrix} x_{2b} \\ x_{3b} \\ x_4 \end{pmatrix} = \begin{pmatrix} 1.2 \\ 1.4 \\ 4 \end{pmatrix} + \begin{pmatrix} \omega_{2b} \\ \omega_{3b} \\ 0 \end{pmatrix}$$

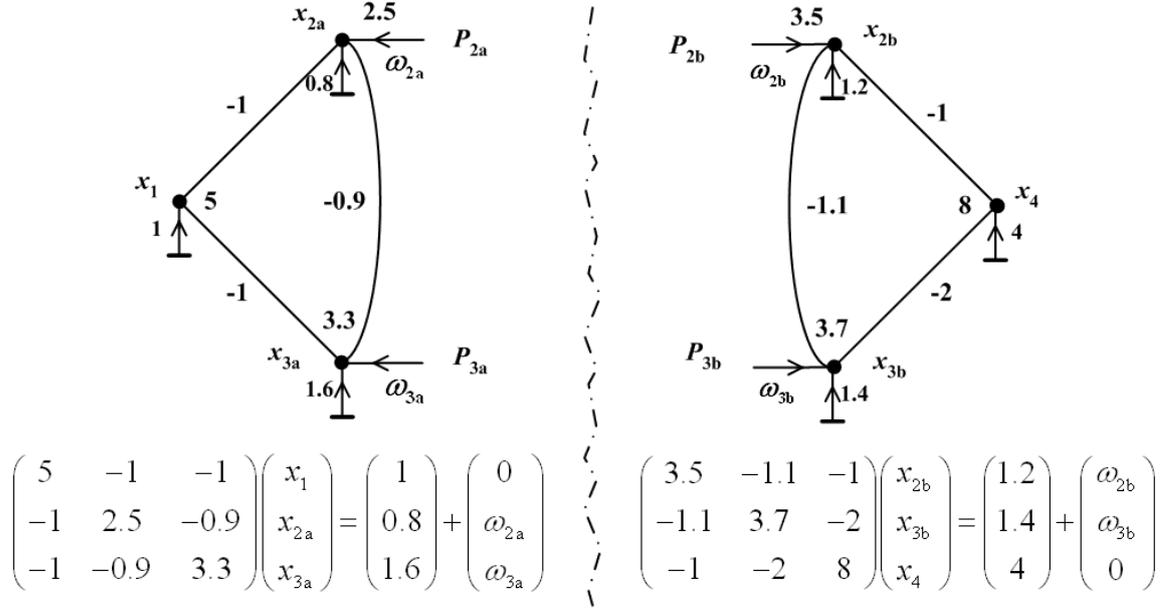

**Figure 5.** Illustration of EVS.

Assume the original graph $G_e$ is partitioned into $N$ separated subgraphs, $M_j, j = 1, 2, \cdots, N$. For two subgraphs, if each of them has at least one port belonging to the same pair of twin vertices, they are called adjacent subgraphs.

According to Section 3, each subgraph could be mapped back into a linear system. To express this linear system, we define $\Gamma_{j,port}$ to be an ordered set of the ports in $M_j$, and $\Gamma_{j,inner}$ to be an ordered set of the inner vertices in $M_j$. Then, we define $\mathbf{u}_j$ to be the potential vector of $\Gamma_{j,port}$, and $\mathbf{y}_j$ to be the potential vector of $\Gamma_{j,inner}$. Then, the subsystem for each subgraph could be expressed by the following equation:

(4.3)
$$\begin{bmatrix} \mathbf{C}_j & \mathbf{E}_j \\ \mathbf{F}_j & \mathbf{D}_j \end{bmatrix} \begin{bmatrix} \mathbf{u}_j \\ \mathbf{y}_j \end{bmatrix} = \begin{bmatrix} \mathbf{f}_j \\ \mathbf{g}_j \end{bmatrix} + \begin{bmatrix} \boldsymbol{\omega}_j \\ 0 \end{bmatrix}$$

where $j = 1, 2, \cdots, N$. $\boldsymbol{\omega}_j$ is the inflow current vector of ports. The inflow current of an inner vertex is zero. $\mathbf{u}_j$ and $\boldsymbol{\omega}_j$ are called the local boundary conditions of $M_j$.

It should be noted that the split vertices could be split again and again, which are called multilevel wire tearing, as seen in Fig. 6, in contrast to the level-one wire tearing mentioned above. The level-two or level-three wire tearing might be applied to partitioning 2- or 3-dimention physical problems. Because of the limited space, we mainly focus on the level-one wire tearing.

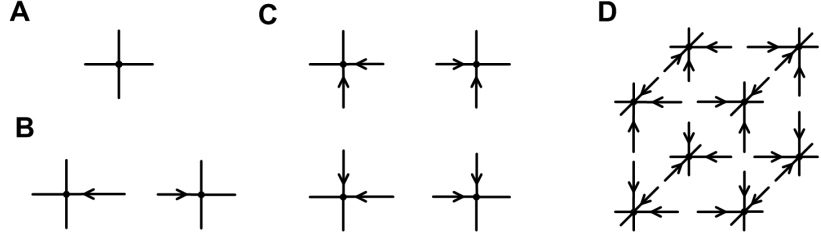

**Figure 6.** Illustration of the multilevel wire tearing. (**A**) The original vertex. (**B**) The twin vertices after the level-one splitting. (**C**) Four vertices after the level-two splitting. (**D**) Eight vertices after the level-three splitting.

## 5. DTM

Assume that the electric graph $G_e$ has been partitioned into $N$ subgraphs, then we insert one DTLP between each pair of twin vertices, which means that we use the directed Transmission Delay Equations as the distributedly-iterative relationship between the boundary conditions of the twin vertices. A simple example is given as below.

**Example 5.1:** This example is based on Example 4.2. After splitting the original graph into two subgraphs, we are going to compute this problem using two processors. Subgraph 1 is located at Processor A, and Subgraph 2 is located at Processor B. Here we assume that the communication delay from Processor A to B is 6.7 μs and that from Processor B to A is 2.9 μs. Fig. 7A illustrates the architecture of this simple parallel computer.

We insert one DTLP between $V_{2a}$ and $V_{2b}$, as shown in Fig. 7A. Their characteristic impedance $Z_2$ between $V_{2a}$ and $V_{2b}$ is set to be 0.2. Then we insert another pair between $V_{3a}$ and $V_{3b}$, and the characteristic impedance $Z_3$ of them are set to be 0.1.

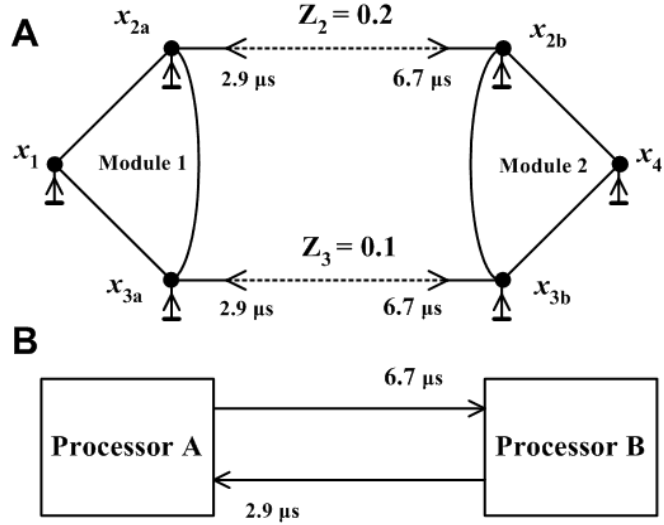

**Figure 7.** Illustration of the algorithm-architecture delay mapping of DTM. (**A**)

The architecture of the parallel computer. **(B)** Illustration of inserting DTLP between each pair of twin vertices.

The propagation delays of DTL from Subgraph 1 to 2 are set to be 6.7 μs, which is the same as the communication delay from Processor A to B, and the propagation delays of DTL from Subgraph 2 to 1 are set to be 2.9 μs. This is an instance of the algorithm-architecture delay mapping, as mentioned in Section 1 and 2.

According to (2.2), the mathematical equation of the DTLP between $V_{2a}$ and $V_{2b}$ is:

(5.1) $$\begin{cases} x_{2a}(t) + 0.2 \cdot \omega_{2a}(t) = x_{2b}(t-2.9) - 0.2 \cdot \omega_{2b}(t-2.9) \\ x_{2b}(t) + 0.2 \cdot \omega_{2b}(t) = x_{2a}(t-6.7) - 0.2 \cdot \omega_{2a}(t-6.7) \end{cases}$$

And the mathematical equation of the DTLP between $V_{3a}$ and $V_{3b}$ is:

(5.2) $$\begin{cases} x_{3a}(t) + 0.1 \cdot \omega_{3a}(t) = x_{3b}(t-2.9) - 0.1 \cdot \omega_{3b}(t-2.9) \\ x_{3b}(t) + 0.1 \cdot \omega_{3b}(t) = x_{3a}(t-6.7) - 0.1 \cdot \omega_{3a}(t-6.7) \end{cases}$$

With (4.1), (5.1) and (5.2), the linear system of Subgraph 1, running on Processor A, could be expressed as below:

(5.3) $$\begin{cases} \begin{pmatrix} 5 & -1 & -1 \\ -1 & 2.5 & -0.9 \\ -1 & -0.9 & 3.3 \end{pmatrix} \begin{pmatrix} x_1(t) \\ x_{2a}(t) \\ x_{3a}(t) \end{pmatrix} = \begin{pmatrix} 1 \\ 0.8 \\ 1.6 \end{pmatrix} + \begin{pmatrix} 0 \\ \omega_{2a}(t) \\ \omega_{3a}(t) \end{pmatrix} \\ x_{2a}(t) + 0.2 \cdot \omega_{2a}(t) = x_{2b}(t-2.9) - 0.2 \cdot \omega_{2b}(t-2.9) \\ x_{3a}(t) + 0.1 \cdot \omega_{3a}(t) = x_{3b}(t-2.9) - 0.1 \cdot \omega_{3b}(t-2.9) \end{cases}$$

Eliminate $\omega_{2a}(t)$ and $\omega_{2b}(t)$ from (5.3), we get following simplified description,

(5.4) $$\begin{pmatrix} 5 & -1 & -1 \\ -1 & 7.5 & -0.9 \\ -1 & -0.9 & 13.3 \end{pmatrix} \begin{pmatrix} x_1(t) \\ x_{2a}(t) \\ x_{3a}(t) \end{pmatrix} = \begin{pmatrix} 1 \\ 0.8 + 5 \cdot x_{2b}(t-2.9) - \omega_{2b}(t-2.9) \\ 1.6 + 10 \cdot x_{3b}(t-2.9) - \omega_{3b}(t-2.9) \end{pmatrix}$$

$$\begin{cases} \omega_{2a}(t) = -5x_{2a}(t) + 5x_{2b}(t-2.9) - \omega_{2b}(t-2.9) \\ \omega_{3a}(t) = -10x_{3a}(t) + 10x_{3b}(t-2.9) - \omega_{3b}(t-2.9) \end{cases}$$

Using the same way, we could combine (4.2), (5.1) and (5.2) to get the simplified description of Subgraph 2 running on Processor B:

(5.5) $$\begin{pmatrix} 8.5 & -1.1 & -1 \\ -1.1 & 13.7 & -2 \\ -1 & -2 & 8 \end{pmatrix} \begin{pmatrix} x_{2b}(t) \\ x_{3b}(t) \\ x_4(t) \end{pmatrix} = \begin{pmatrix} 1.2 + 5 \cdot x_{2a}(t-6.7) - \omega_{2a}(t-6.7) \\ 1.4 + 10 \cdot x_{3a}(t-6.7) - \omega_{3a}(t-6.7) \\ 4 \end{pmatrix}$$

$$\begin{cases} \omega_{2b}(t) = -5x_{2b}(t) + 5x_{2a}(t-6.7) - \omega_{2a}(t-6.7) \\ \omega_{3b}(t) = -10x_{3b}(t) + 10x_{3a}(t-6.7) - \omega_{3a}(t-6.7) \end{cases}$$

At the right hand side of (5.4), the time-delay variables are the previous computing result received from Processor B, so they are known at the current time. Processor A solves (5.4) to get the potentials and currents of Subgraph 1 at the current time, and then send them to Processor B. This is the computing process of Processor A, and so is that of Processor B.

At last, we set the initial value of (5.4) and (5.5) in (5.6) and do the asynchronous and distributed computing. The result is shown in Fig. 8.

$$(5.6) \quad \begin{cases} x_{2a}(0) = x_{2b}(0) = x_{3a}(0) = x_{3b}(0) = 0 \\ \omega_{2a}(0) = \omega_{2b}(0) = \omega_{3a}(0) = \omega_{3b}(0) = 0 \end{cases}$$

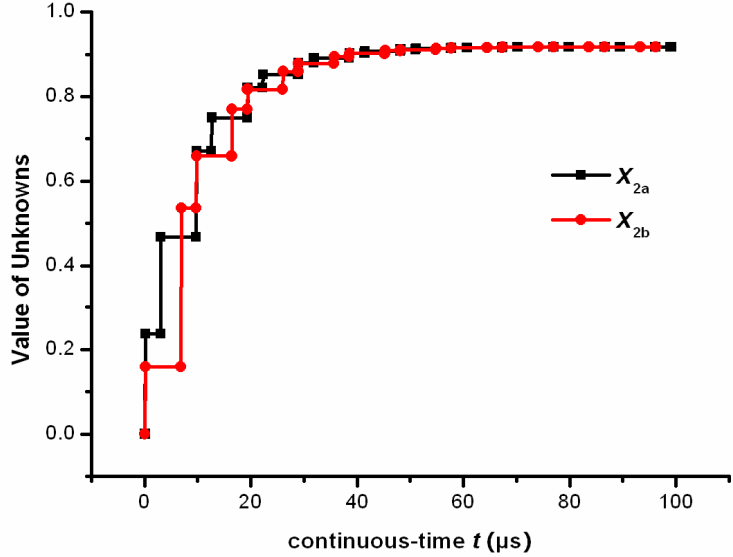

**Figure 8.** Computing result of DTM

The choice of the characteristic impedance of DTLP, i.e. $Z_2$ and $Z_3$, would affect the convergence speed of the algorithm. Fig. 9 illustrates this effect. As the result, we could speedup DTM if the characteristic impedances of DTLPs are carefully chosen.

After illustrating DTM by an example, we give the mathematical description of this algorithm. Assume that the electric graph $G_e$ has been partitioned into $N$ subgraphs, $M_j, j = 1, 2, \cdots, N$, then we insert one DTLP between each pair of twin vertices. This means that we use the Directed Transmission Delay Equations as the boundary conditions for (4.3).

In the previous section, we have defined $\Gamma_{j,port}$ as an ordered set of the ports of $M_j$, and we defined $\mathbf{u}_j$ to be the potential vector of $\Gamma_{j,port}$, and $\mathbf{y}_j$ to be the potential vector of $\Gamma_{j,inner}$. Further, we define $\Gamma_{j,twin}$ to be another ordered set of ports whose twin vertices belong to $\Gamma_{j,port}$. The vertices in $\Gamma_{j,port}$ and their corresponding twin vertices in $\Gamma_{j,twin}$ have the same order. The vertices of $\Gamma_{j,twin}$ belong to the adjacent subgraphs of $M_j$.

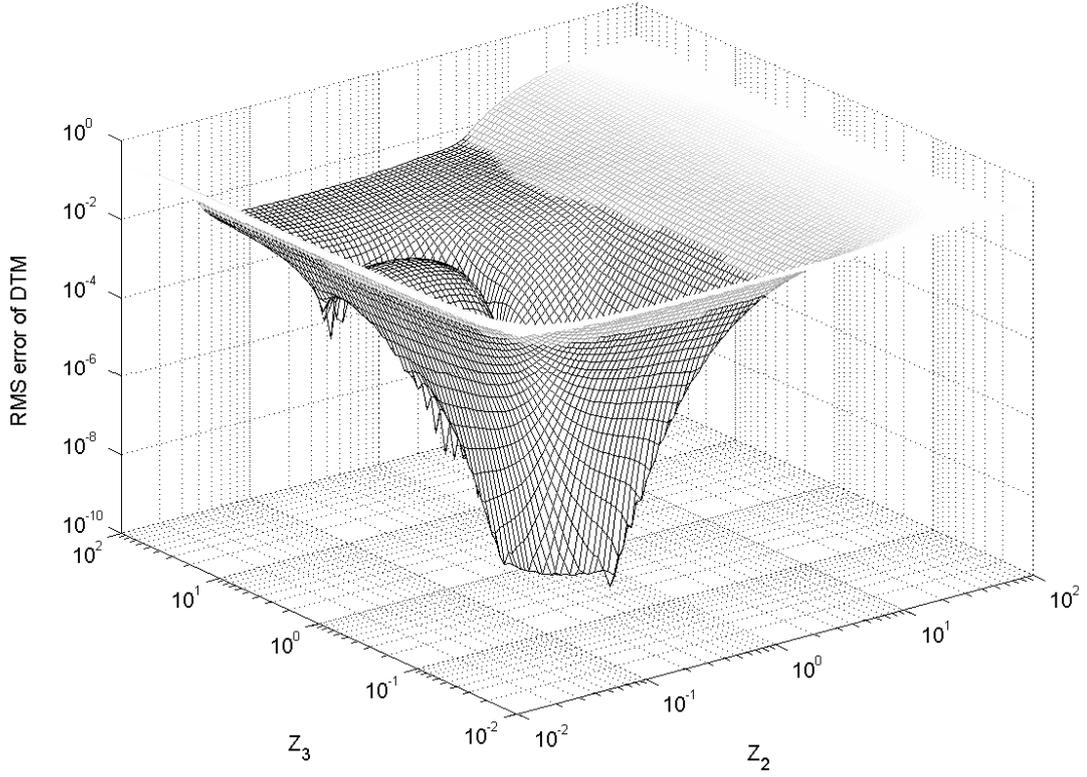

**Figure 9.** RMS error of DTM when $t = 100$ μs

Define $\mathbf{u}_{j,twin}$ as the potential vector of $\Gamma_{j,twin}$, and $\mathbf{\omega}_{j,twin}$ as the current vector of $\Gamma_{j,twin}$. Then, (5.7) expresses the Directed Transmission Delay Equations of DTL pointed to $M_j$. Note that DTL apart from $M_j$ are not considered here, and they would be considered by the subgraphs they point to.

(5.7) $\quad\quad\quad\quad \mathbf{u}_j(t) + \mathbf{Z}_j \cdot \mathbf{\omega}_j(t) = \mathbf{u}_{j,twin}(t-\boldsymbol{\tau}) - \mathbf{Z}_j \cdot \mathbf{\omega}_{j,twin}(t-\boldsymbol{\tau})$

Where $t$ is the continuous-time variable. $r$ is the total number of DTL pointed to $M_j$. $\tau_i$, $i = 1,\cdots,r$, is the propagation delays of the $i$-th DTL.

$$\mathbf{u}_j(t) = \begin{bmatrix} u_1(t) & u_2(t) & \cdots & u_r(t) \end{bmatrix}^{\mathrm{T}},$$

$$\mathbf{\omega}_j(t) = \begin{bmatrix} \omega_1(t) & \omega_2(t) & \cdots & \omega_r(t) \end{bmatrix}^{\mathrm{T}},$$

$$\mathbf{u}_{j,twin}(t-\boldsymbol{\tau}) = \begin{bmatrix} u_{1,twin}(t-\tau_1) & u_{2,twin}(t-\tau_2) & \cdots & u_{r,twin}(t-\tau_r) \end{bmatrix}^{\mathrm{T}},$$

$$\mathbf{\omega}_{j,twin}(t-\boldsymbol{\tau}) = \begin{bmatrix} \omega_{1,twin}(t-\tau_1) & \omega_{2,twin}(t-\tau_2) & \cdots & \omega_{r,twin}(t-\tau_r) \end{bmatrix}^{\mathrm{T}}.$$

$\mathbf{Z}_j = diag(z_1,\cdots,z_r)$, is a positive diagonal matrix, called the local characteristic impedance matrix of $M_j$. Its diagonal elements are the characteristic impedances of DTL. Note that the characteristic impedances of the DTL belonging to the same DTLP should be same.

(5.7) is a continuous-time iterative relation, and $\mathbf{u}_{j,twin}(t-\tau)$ and $\boldsymbol{\omega}_{j,twin}(t-\tau)$ are the previous computing results passed from the adjacent subgraphs of $M_j$, which are called the remote boundary condition of $M_j$. Merge (5.7) and (4.3), we get:

$$(5.8) \quad \begin{bmatrix} \mathbf{C}_j & \mathbf{E}_j & -\mathbf{I} \\ \mathbf{F}_j & \mathbf{D}_j & 0 \\ \mathbf{I} & 0 & \mathbf{Z}_j \end{bmatrix} \begin{bmatrix} \mathbf{u}_j(t) \\ \mathbf{y}_j(t) \\ \boldsymbol{\omega}_j(t) \end{bmatrix} = \begin{bmatrix} \mathbf{f}_j \\ \mathbf{g}_j \\ \mathbf{u}_{j,twin}(t-\tau) - \mathbf{Z}_j \cdot \boldsymbol{\omega}_{j,twin}(t-\tau) \end{bmatrix}$$

where $\mathbf{I}$ is the identity matrix.

Removing $\boldsymbol{\omega}_j(t)$ from (5.8), we obtain the following SPD system.

$$(5.9) \quad \begin{bmatrix} \mathbf{C}_j + \mathbf{Z}_j^{-1} & \mathbf{E}_j \\ \mathbf{F}_j & \mathbf{D}_j \end{bmatrix} \begin{bmatrix} \mathbf{u}_j(t) \\ \mathbf{y}_j(t) \end{bmatrix} = \begin{bmatrix} \mathbf{f}_j + \mathbf{Z}_j^{-1} \{\mathbf{u}_{j,twin}(t-\tau) - \mathbf{Z}_j \cdot \boldsymbol{\omega}_{j,twin}(t-\tau)\} \\ \mathbf{g}_j \end{bmatrix}$$

$$\boldsymbol{\omega}_j(t) = -\mathbf{Z}^{-1}\mathbf{u}_j(t) + \mathbf{Z}^{-1}\mathbf{u}_{j,twin}(t-\tau) - \boldsymbol{\omega}_{j,twin}(t-\tau)$$

(5.9) is called the local system of $M_j$, which should be solved when the remote boundary condition is updated. (5.9) could be solved by Sparse or Dense Cholesky, CG, MG, etc.

It is noticeable that the coefficient matrix of (5.9) is constant during computing process, and this is the key to speed up DTM. For instance, if we use the Cholesky factorization to solve the local system, actually only once factorization should be done at the beginning; as long as we get the Cholesky factor, it is a piece of cake to solve (5.9) since we just need to do the forward and backward substitution in the following time.

If all the DTLs have an equal propagation delay, noted by 1 time unit, then (5.9) could be re-expressed as the discrete-time iterative form (5.10), which is the local system of VTM [6].

$$(5.10) \quad \begin{bmatrix} \mathbf{C}_j + \mathbf{Z}_j^{-1} & \mathbf{E}_j \\ \mathbf{F}_j & \mathbf{D}_j \end{bmatrix} \begin{bmatrix} \mathbf{u}_j^k \\ \mathbf{y}_j^k \end{bmatrix} = \begin{bmatrix} \mathbf{f}_j + \mathbf{Z}_j^{-1}(\mathbf{u}_{j,twin}^{k-1} - \mathbf{Z}_j \boldsymbol{\omega}_{j,twin}^{k-1}) \\ \mathbf{g}_j \end{bmatrix}$$

$$\boldsymbol{\omega}_j^k = -\mathbf{Z}^{-1} \cdot \mathbf{u}_j^k + \mathbf{Z}^{-1} \cdot \mathbf{u}_{j,twin}^{k-1} - \boldsymbol{\omega}_{j,twin}^{k-1}$$

Referring to (5.4), (5.5) and (5.9), it is noticeable that DTM is different from the traditional iterative algorithms, i.e. Gauss Jacobi, which usually have the discrete-time iterative form, i.e. $\mathbf{x}^k = f(\mathbf{x}^{k-1})$. In DTM, the continuous-time variable $t$ is used instead of the iterative index $k$. This indicates that DTM is a continuous-time iterative algorithm, which is more flexible than the traditional discrete-time algorithms.

Table 1 gives the detail of the DTM algorithm. It should be noted that there is no synchronization step, no broadcasting, but only N2N communication. Once one processor receives the remote boundary conditions from one or more adjacent processors, it could immediately do its local computation without waiting for the slowest processor. DTM is suited to be implemented using the message passing approach [13].

**Table 1. Algorithm of DTM**

Assume the original electric graph is partitioned into $N$ subgraphs. For Subgraph $M_j$, $j = 1, \cdots, N$, do in parallel:

1. Guess the initial local boundary condition, $\mathbf{u}_j(0)$ and $\mathbf{\omega}_j(0)$, of each port.

2. Communicate with adjacent subgraphs, to make an agreement of the characteristic impedances for each DTLP. As the result, $\mathbf{Z}_j$ is set.

3. Wait until receiving part of the remote boundary conditions, $\mathbf{u}_{j,twin}(t-\tau)$ and $\mathbf{\omega}_{j,twin}(t-\tau)$, from one or more of the adjacent subgraphs.

    3.1 Solve the local system with the updated remote boundary condition and obtain the new local boundary condition, $\mathbf{u}_j(t)$ and $\mathbf{\omega}_j(t)$.

    3.2 Send the new local boundary condition to the adjacent subgraphs of $M_j$.

    3.3. If convergent, then break.

4. EndWait.

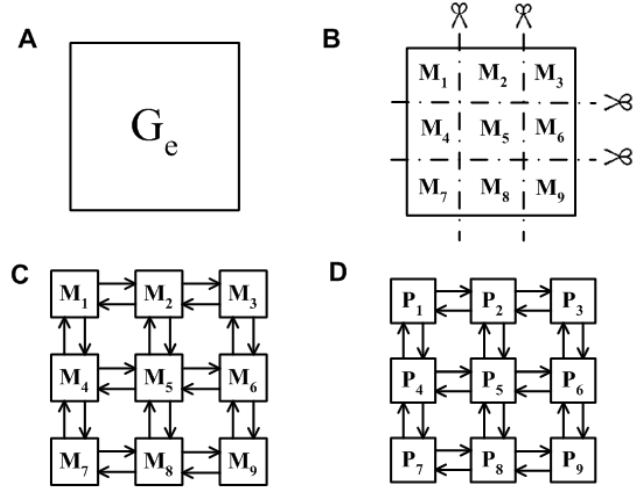

**Figure 10.** Illustration of the computing process of DTM. (**A**) The original electric graph of the sparse linear system. (**B**) Partition the original graph into $N$ subgraphs by EVS. (**C**) Add DTLs between adjacent subgraphs. (**D**) Map each subgraph onto one processor, and map each DTL to a directed communication path.

## 6. Convergence theory

**Theorem 6.1** (Convergence): Assume an SPD linear system is partitioned into $N$ subgraphs by EVS. If there is at least one SPD subgraph, and the other subgraphs are symmetric-non-negative-definite (SNND), then DTM converges to the solution of the original system. The characteristic impedances of DTLP could be set to arbitrary positive values.

This conclusion is valid for both the level-one and multilevel wire tearing. For the case of level-one splitting, a basic proof is given in the Appendix according to the final-value theorem of Laplacian transformation [14].

## 7. Numerical experiments

In this section, we test DTM on DTM toolbox, which is developed by us under MATLAB & SIMULINK [15, 16]. This toolbox gives us an easy way to simulate the behavior of DTM on heterogeneous parallel platform.

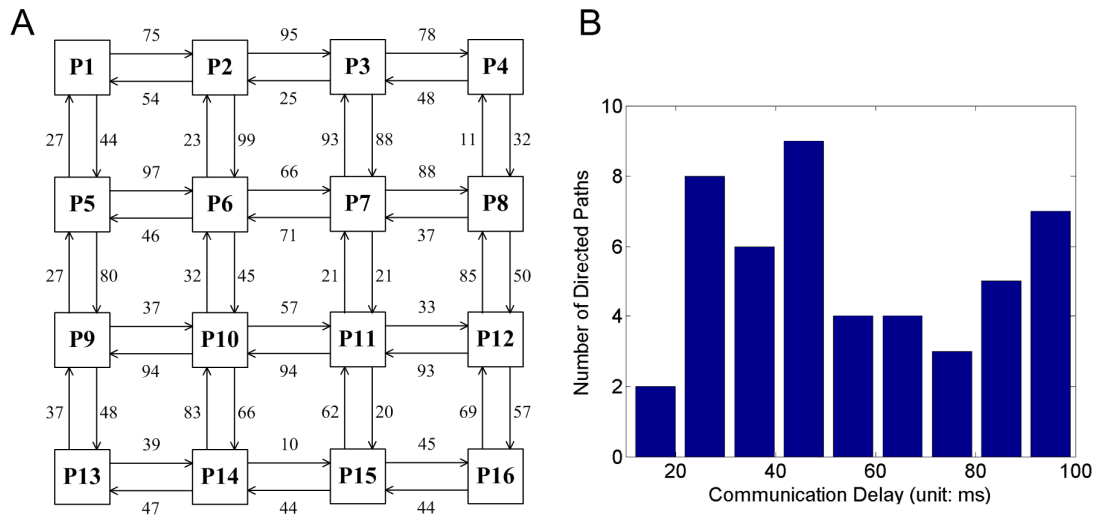

**Figure 11.** (**A**) Heterogeneous topology of 16 processors in a 4x4 mesh. The N2N communication delays are also illustrated (unit: ms). (**B**) Bar chart of the N2N communication delays.

We first test DTM on 16 processors configured as a 4x4 mesh shown in Fig. 11. Here the communication network is very unsymmetrical. The maximum delay (99ms) is about 9 times larger than the minimum delay (10ms) and the delay from Processor $P_k$ to $P_j$ is quite different from the delay from Processor $P_j$ to $P_k$.

This is a terrible parallel environment for the parallel algorithms which need synchronization and broadcasting, while it is cozy for DTM. Fig. 12 shows the convergence curve of DTM to solve linear systems on 16 processors. The sparse SPD linear systems for test are randomly generated and regularly partitioned using the level-one and level-two mixed EVS to achieve load balance.

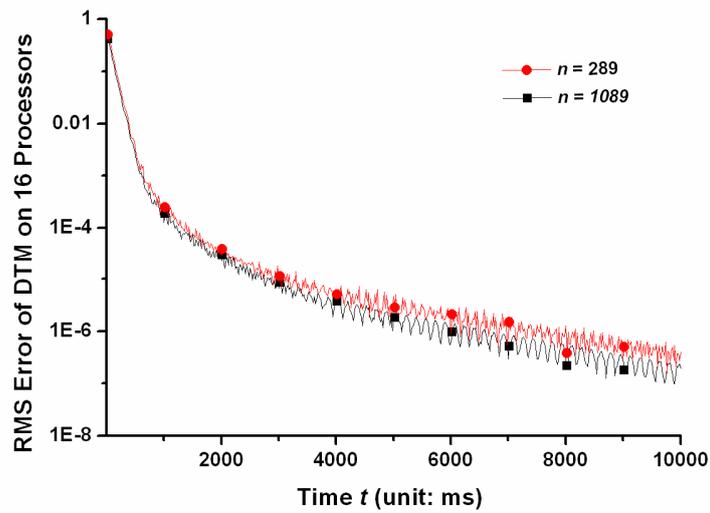

**Figure 12.** Computing result of DTM on 16 processors

Then we test DTM on 64 processors configured as an 8x8 mesh, as shown in Fig. 13. The N2N communication delays are uniformly distributed between 10ms and 100ms. Fig. 14 illustrates the computational errors of two sparse linear systems, having 1089 and 4225 unknowns, respectively.

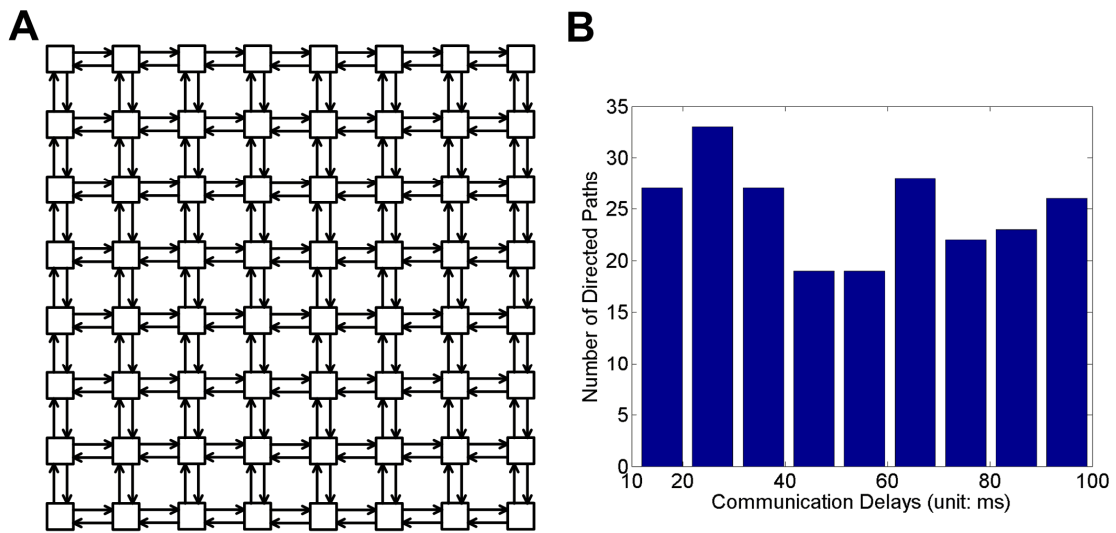

**Figure 13.** (**A**) Heterogeneous topology of the 64 processors in an 8×8 mesh. The N2N communication delays are different. (**B**) Bar chart of the N2N communication delays.

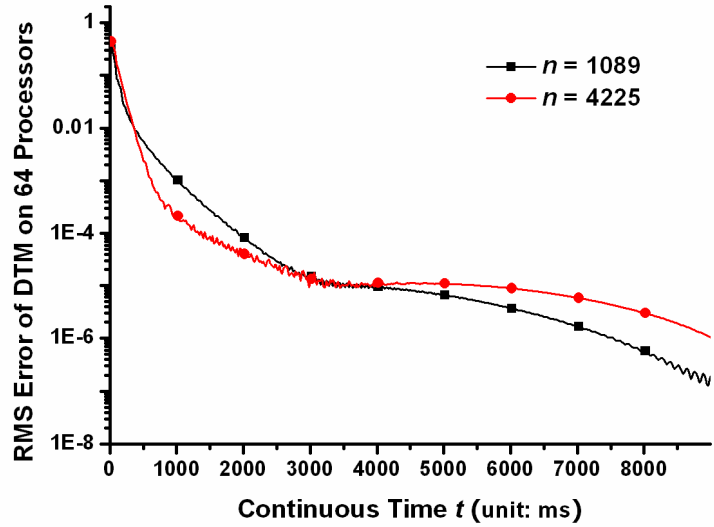

**Figure 14.** Computing result of DTM on 64 processors.

Theoretically, the dimension of the sparse linear system being solved by DTM could be arbitrarily-large, and the processors could be arbitrary number. Limited by our hardware, here we only test *n* = 289, 1089 and 4225 on 16 and 64 processors.

## 8. Conclusions

In this paper, we propose a new parallel algorithm, DTM, to solve the sparse SPD linear systems, and we bring in EVS to partition the electric graph of the symmetric linear system. We present the convergence theorem, which makes DTM feasible for any kinds of SPD systems.

DTM is an asynchronous, distributed and continuous-time iterative algorithm and it is able to be freely running on all kinds of heterogeneous parallel computers, e.g. multicore or manycore microprocessor, clusters, grids, clouds, Internet, ad hoc network and wireless network.

Compared to VTM, the convergence speed of DTM is slower. We wonder if there is some way to lessen the speed gap between DTM and VTM. Since VTM could be considered as a global synchronous version of DTM, we guess that a tradeoff between DTM and VTM could be made by some sync-async-mixed approach in the physical domain (e.g. global-async-local-sync) or time domain (e.g. async-sync-async-sync, synchronizing once after a period of asynchronization).

## Acknowledgments

We thank Prof. Hao Zhang, Yi Su, Dr. Chun Xia, Dr. Wei Xue, Pei Yang, Bin Niu. This work was partially sponsored by the Major State Basic Research Development Program of China (973 Program) under contract G1999032903, the National Natural Science Foundation of China Key Program under Grant #90207001, and the National Science Fund for Distinguished Young Scholars of China under Grant #60025101.

# References

[1]. J. W. Demmel. Lecture notes of Application of Parallel Computers at University of California, Berkeley. Available at http://www.cs.berkeley.edu/~demmel/cs267/

[2]. A. George and J. W. Liu. Computer solution of sparse systems, Prentice Hall, 1981.

[3]. Y. Saad. Iterative methods for sparse linear systems, 2$^{nd}$ edition, SIAM, 2003.

[4]. R. Barrett, M. Berry, T. Chan, et al. Templates for the solution of Linear Systems: Building Blocks for Iterative Methods, 2$^{nd}$ Edition, SIAM, 1994.

[5]. A. Toselli and O. Widlund. Domain decomposition methods - algorithms and theory, Springer, 2005.

[6]. F. Wei and H. Yang. Virtual Transmission Method, a new distributed algorithm to solve sparse linear systems. NCM, 2008.

[7]. M. B. Taylor, J. Kim, J. Miller, et al. The RAW microprocessor: a computational fabric for software circuits and general-purpose program. IEEE Micro. 22, 2 (March, 2002), 25-35.

[8]. S. Dighe, H. Wilson, J. Tschanz, et al. An 80-Tile 1.28 TFLOPS Network-on-Chip in 65nm CMOS. In Proceedings of International Solid-State Circuits Conference, 2007.

[9]. H. J. Pain. The physics of vibrations and waves, Wiley, 1976.

[10]. R. E. Collin. Foundations for microwave engineering, 2$^{nd}$ edition, Wiley-IEEE Press, 2000.

[11]. J. H. Gridley. Principles of electrical transmission lines in power and communications, Pergamon Press, 1967.

[12]. T. L. Floyd. Principles of electric circuits, 6$^{th}$ edition, Prentice Hall, 1999.

[13]. A. Grama, A. Gupta, G. Karypis and V. Kumar, Introduction to parallel computing, 2$^{nd}$ edition, Addison-Wesley, 2002.

[14]. A. V. Oppenheim, A. S. Willsky and S. H. Nawab. Signals and Systems. 2$^{nd}$ edition. Englewood Cliffs, NJ: Prentice Hall. 1996.

[15] Matlab User Manual, Version 6.5 R13, The Math Works Inc., 2002.

[16] Simulink References, Version 5.0 R13, The Math Works Inc., 2002.

[17] Gérard M. Baudet, Asynchronous Iterative Methods for Multiprocessors, Journal of the ACM (JACM), April 1978.

[18] Andreas Frommer, Daniel B. Szyld, On asynchronous iterations, Journal of Computational and Applied Mathematics, 2000.

[19] J. Bahi, S. Contassot-Vivier, R. Couturier. Parallel Iterative Algorithms: From Sequential to Grid Computing. Chapman & Hall/CRC, 2007.

[20] Fei Wei, Huazhong Yang. Directed Transmission Method, a fully asynchronous approach to solve sparse linear system in parallel, poster, ACM SPAA, June 2008.

[16] Felix F.Wu. Solution of large-scale networks by tearing. IEEE Transactions on Circuits and Systems, 1976.

[17] Alberto Sangiovanni-Vincentelli, Li-Kuan Chen, and Leon O. Chua. A new tearing approach – node-tearing nodal analysis. In IEEE International Symposium on Circuits and Systems, 1977.

[18] P. Yang, I. N. Hajj, and T. N. Trick. SLATE: A circuit simulation program with latency exploitation and node tearing. In IEEE International Conference on Circuits and Computers, October 1980.

[19] N. B. Guy Rabbat, Alberto L. Sangiovanni-Vincentelli, and Hsueh Y. Hsieh. A multilevel newton algorithm with macromodeling and latency for the analysis of large-scale nonlinear circuits in the time domain. IEEE Transactions on Circuits and Systems, 26(9):733–741, September 1979.

## Appendix.   Proof of the convergence theorem

Here we prove a simple version of the convergence theory, and we assume that there is no inner vertex and all the vertices are split. Actually, if there were inner vertices, they could be eliminated and we still able to get (A.1).

(A.1) $$\mathbf{A} \cdot \mathbf{u} = \mathbf{b}$$

Assume that the SPD system (A.1) is partitioned into two SPD subgraphs (A.2) by level-one wire tearing technique, and the characteristic impedances of DTLPs are set to be arbitrary positive value.

(A.2) $$\begin{cases} \mathbf{A}_1 \mathbf{u}_1(t) = \mathbf{b}_1 + \boldsymbol{\omega}_1(t) \\ \mathbf{A}_2 \mathbf{u}_2(t) = \mathbf{b}_2 + \boldsymbol{\omega}_2(t) \end{cases}$$

where $\mathbf{A} = \mathbf{A}_1 + \mathbf{A}_2$, $\mathbf{b} = \mathbf{b}_1 + \mathbf{b}_2$. The DTLs could be described as below:

(A.3) $$\begin{cases} \mathbf{u}_1(t) + \mathbf{Z} \cdot \boldsymbol{\omega}_1(t) = \mathbf{u}_2(t-\boldsymbol{\tau}) - \mathbf{Z} \cdot \boldsymbol{\omega}_2(t-\boldsymbol{\tau}) \\ \mathbf{u}_2(t) + \mathbf{Z} \cdot \boldsymbol{\omega}_2(t) = \mathbf{u}_1(t-\boldsymbol{\sigma}) - \mathbf{Z} \cdot \boldsymbol{\omega}_1(t-\boldsymbol{\sigma}) \end{cases}$$

where $\mathbf{Z} = diag(z_1, \cdots, z_r)$, which is the characteristic impedance matrix. Below is another description of convergence theorem when $N = 2$.

**Lemma A.1**: Suppose the solutions of (A.2) by DTM on two parallel processors are $\mathbf{u}_1(t)$ and $\mathbf{u}_2(t)$. $\forall z_i > 0, i = 1, 2, \cdots, r$, we have

$$\lim_{t \to +\infty} \mathbf{u}_1(t) = \lim_{t \to +\infty} \mathbf{u}_2(t) = \mathbf{A}^{-1}\mathbf{b}$$

To prove Lemma A.1, we prove following Lemma at first.

**Lemma A.2**: $\mathbf{A}$ is a $r \times r$ SPD martrix, $\mathbf{Z}$ is a $r \times r$ positive diagonal martrix, $\exists \mathbf{Q}$, $\mathbf{Q}\mathbf{Q}^T = \mathbf{I}$, and $\mathbf{ZA} = \sqrt{\mathbf{Z}}\mathbf{Q}\mathbf{T}\mathbf{Q}^T\sqrt{\mathbf{Z}}^{-1}$, where $\mathbf{T} = diag(t_1, t_2, \cdots t_r)$, $t_i$ is the $i$-th eigen value of $\mathbf{ZA}$.

**Proof:**

As $\sqrt{\mathbf{Z}}\mathbf{A}\sqrt{\mathbf{Z}}$ is symmetric and positive definite, $\exists \mathbf{Q}$, $\mathbf{Q}\mathbf{Q}^T = \mathbf{I}$, and $\sqrt{\mathbf{Z}}\mathbf{A}\sqrt{\mathbf{Z}} = \mathbf{Q}\mathbf{T}\mathbf{Q}^T$, where $\mathbf{T} = diag(t_1, t_2, \cdots t_r)$, $t_i$ is the $i$-th eigen value of $\sqrt{\mathbf{Z}}\mathbf{A}\sqrt{\mathbf{Z}}$. Obviously, $\sqrt{\mathbf{Z}}\mathbf{A}\sqrt{\mathbf{Z}} = \mathbf{ZA}$, thus $t_i$ is the $i$-th eigenvalue of $\mathbf{ZA}$ too. Therefore we have

$$\mathbf{ZA} = \sqrt{\mathbf{Z}}\sqrt{\mathbf{Z}}\mathbf{A}\sqrt{\mathbf{Z}}\sqrt{\mathbf{Z}}^{-1} = \sqrt{\mathbf{Z}}\mathbf{Q}\mathbf{T}\mathbf{Q}^T\sqrt{\mathbf{Z}}^{-1}$$

∎

**Proof of Lemma A.1:**

We perform Laplacian transformation to (A.2) at first

(A.4)
$$\begin{cases} \mathbf{A}_1 \cdot \mathbf{U}_1(s) = \dfrac{\mathbf{b}_1}{s} + \mathbf{\Omega}_1(s) \\ \mathbf{A}_2 \cdot \mathbf{U}_2(s) = \dfrac{\mathbf{b}_2}{s} + \mathbf{\Omega}_2(s) \end{cases}$$

(A.5)
$$\begin{cases} \mathbf{U}_1(s) + \mathbf{Z} \cdot \mathbf{\Omega}_1(s) = e^{-s\boldsymbol{\tau}} \cdot \mathbf{U}_2(s) - \mathbf{Z} \cdot e^{-s\boldsymbol{\tau}} \cdot \mathbf{\Omega}_2(s) \\ \mathbf{U}_2(s) + \mathbf{Z} \cdot \mathbf{\Omega}_2(s) = e^{-s\boldsymbol{\sigma}} \cdot \mathbf{U}_1(s) - \mathbf{Z} \cdot e^{-s\boldsymbol{\sigma}} \cdot \mathbf{\Omega}_1(s) \end{cases}$$

Where $e^{-s\boldsymbol{\tau}} = diag(e^{-s\tau_1}, \cdots, e^{-s\tau_r})$. $e^{-s\boldsymbol{\sigma}} = diag(e^{-s\sigma_1}, \cdots, e^{-s\sigma_r})$. $r$ is the total number of DTLP.

Remove $\mathbf{\Omega}_1(s)$ and $\mathbf{\Omega}_2(s)$ from (A.5), and we get:

$$\begin{cases} (\mathbf{I} + \mathbf{Z} \cdot \mathbf{A}_1) \cdot \mathbf{U}_1(s) - \dfrac{\mathbf{Z} \cdot \mathbf{b}_1}{s} = e^{-s\boldsymbol{\tau}} \cdot (\mathbf{I} - \mathbf{Z} \cdot \mathbf{A}_2) \cdot \mathbf{U}_2(s) + e^{-s\boldsymbol{\tau}} \cdot \dfrac{\mathbf{Z} \cdot \mathbf{b}_2}{s} \\ (\mathbf{I} + \mathbf{Z} \cdot \mathbf{A}_2) \cdot \mathbf{U}_2(s) - \dfrac{\mathbf{Z} \cdot \mathbf{b}_2}{s} = e^{-s\boldsymbol{\sigma}} \cdot (\mathbf{I} - \mathbf{Z} \cdot \mathbf{A}_1) \cdot \mathbf{U}_1(s) + e^{-s\boldsymbol{\sigma}} \cdot \dfrac{\mathbf{Z} \cdot \mathbf{b}_1}{s} \end{cases}$$

$$\begin{cases} \mathbf{U}_1(s) = (\mathbf{I} + \mathbf{Z} \cdot \mathbf{A}_1)^{-1} \cdot e^{-s\boldsymbol{\tau}} \cdot (\mathbf{I} - \mathbf{Z} \cdot \mathbf{A}_2) \cdot \mathbf{U}_2(s) + (\mathbf{I} + \mathbf{Z} \cdot \mathbf{A}_1)^{-1} \cdot \left( e^{-s\boldsymbol{\tau}} \cdot \dfrac{\mathbf{Z} \cdot \mathbf{b}_2}{s} + \dfrac{\mathbf{Z} \cdot \mathbf{b}_1}{s} \right) \\ \mathbf{U}_2(s) = (\mathbf{I} + \mathbf{Z} \cdot \mathbf{A}_2)^{-1} \cdot e^{-s\boldsymbol{\sigma}} \cdot (\mathbf{I} - \mathbf{Z} \cdot \mathbf{A}_1) \cdot \mathbf{U}_1(s) + (\mathbf{I} + \mathbf{Z} \cdot \mathbf{A}_2)^{-1} \cdot \left( e^{-s\boldsymbol{\sigma}} \cdot \dfrac{\mathbf{Z} \cdot \mathbf{b}_1}{s} + \dfrac{\mathbf{Z} \cdot \mathbf{b}_2}{s} \right) \end{cases}$$

Then we get:

(A.6) $\quad \mathbf{U}_1(s) = (\mathbf{I} - \mathbf{T})^{-1} (\mathbf{W}_1 + \mathbf{W}_2)$

Where
$$\mathbf{T} = (\mathbf{I} + \mathbf{Z} \cdot \mathbf{A}_1)^{-1} \cdot e^{-s\boldsymbol{\tau}} \cdot (\mathbf{I} - \mathbf{Z} \cdot \mathbf{A}_2) \cdot (\mathbf{I} + \mathbf{Z} \cdot \mathbf{A}_2)^{-1} \cdot e^{-s\boldsymbol{\sigma}} \cdot (\mathbf{I} - \mathbf{Z} \cdot \mathbf{A}_1)$$

$$\mathbf{W}_1 = (\mathbf{I} + \mathbf{Z} \cdot \mathbf{A}_1)^{-1} \cdot e^{-s\boldsymbol{\tau}} \cdot (\mathbf{I} - \mathbf{Z} \cdot \mathbf{A}_2) \cdot (\mathbf{I} + \mathbf{Z} \cdot \mathbf{A}_2)^{-1} \cdot e^{-s\boldsymbol{\sigma}} \cdot \dfrac{\mathbf{Z} \cdot \mathbf{b}_1}{s}$$
$$+ (\mathbf{I} + \mathbf{Z} \cdot \mathbf{A}_1)^{-1} \dfrac{\mathbf{Z} \cdot \mathbf{b}_1}{s}$$

$$\mathbf{W}_2(s) = (\mathbf{I} + \mathbf{Z} \cdot \mathbf{A}_1)^{-1} \cdot e^{-s\boldsymbol{\tau}} \cdot (\mathbf{I} - \mathbf{Z} \cdot \mathbf{A}_2) \cdot (\mathbf{I} + \mathbf{Z} \cdot \mathbf{A}_2)^{-1} \dfrac{\mathbf{Z} \cdot \mathbf{b}_2}{s}$$
$$+ (\mathbf{I} + \mathbf{Z} \cdot \mathbf{A}_1)^{-1} \cdot e^{-s\boldsymbol{\tau}} \cdot \dfrac{\mathbf{Z} \cdot \mathbf{b}_2}{s}$$

$$\mathbf{I} - \mathbf{T} = \mathbf{I} - (\mathbf{I} + \mathbf{Z} \cdot \mathbf{A}_1)^{-1} \cdot e^{-s\boldsymbol{\tau}} \cdot (\mathbf{I} - \mathbf{Z} \cdot \mathbf{A}_2) \cdot (\mathbf{I} + \mathbf{Z} \cdot \mathbf{A}_2)^{-1} \cdot e^{-s\boldsymbol{\sigma}} \cdot (\mathbf{I} - \mathbf{Z} \cdot \mathbf{A}_1)$$
$$= (\mathbf{I} + \mathbf{Z} \cdot \mathbf{A}_1)^{-1} \mathbf{W} (\mathbf{I} - \mathbf{Z} \cdot \mathbf{A}_1)$$

Where
$$\mathbf{W} = (\mathbf{I} + \mathbf{Z} \cdot \mathbf{A}_1)(\mathbf{I} - \mathbf{Z} \cdot \mathbf{A}_1) - e^{-s\boldsymbol{\tau}} \cdot (\mathbf{I} - \mathbf{Z} \cdot \mathbf{A}_2) \cdot (\mathbf{I} + \mathbf{Z} \cdot \mathbf{A}_2)^{-1} \cdot e^{-s\boldsymbol{\sigma}}.$$

As $\mathbf{A}_1$ and $\mathbf{A}_2$ are SPD, from Lemma A.2 we have,

$$\mathbf{Z} \cdot \mathbf{A}_1 = \sqrt{\mathbf{Z}} \mathbf{Q}_1 \mathbf{T}_1 \mathbf{Q}_1^T \sqrt{\mathbf{Z}}^{-1}$$
$$\mathbf{Z} \cdot \mathbf{A}_2 = \sqrt{\mathbf{Z}} \mathbf{Q}_2 \mathbf{T}_2 \mathbf{Q}_2^T \sqrt{\mathbf{Z}}^{-1}$$

Where $\mathbf{Q}_1 \mathbf{Q}_1^T = \mathbf{Q}_2 \mathbf{Q}_2^T = \mathbf{I}$, $\mathbf{T}_1$ and $\mathbf{T}_2$ and are positive diagonal matrices. Thus we have,

$$\begin{aligned}
\mathbf{W} &= (\mathbf{I}+\mathbf{Z}\cdot\mathbf{A}_1)(\mathbf{I}-\mathbf{Z}\cdot\mathbf{A}_1)^{-1} - e^{-s\boldsymbol{\tau}}\cdot(\mathbf{I}-\mathbf{Z}\cdot\mathbf{A}_2)\cdot(\mathbf{I}+\mathbf{Z}\cdot\mathbf{A}_2)^{-1}\cdot e^{-s\boldsymbol{\sigma}}\\
&= \left(\mathbf{I}+\sqrt{\mathbf{Z}}\mathbf{Q}_1\mathbf{T}_1\mathbf{Q}_1^{\mathrm{T}}\sqrt{\mathbf{Z}}^{-1}\right)\left(\mathbf{I}-\sqrt{\mathbf{Z}}\mathbf{Q}_1\mathbf{T}_1\mathbf{Q}_1^{\mathrm{T}}\sqrt{\mathbf{Z}}^{-1}\right)^{-1}\\
&\quad - e^{-s\boldsymbol{\tau}}\cdot\left(\mathbf{I}-\sqrt{\mathbf{Z}}\mathbf{Q}_2\mathbf{T}_2\mathbf{Q}_2^{\mathrm{T}}\sqrt{\mathbf{Z}}^{-1}\right)\cdot\left(\mathbf{I}+\sqrt{\mathbf{Z}}\mathbf{Q}_2\mathbf{T}_2\mathbf{Q}_2^{\mathrm{T}}\sqrt{\mathbf{Z}}^{-1}\right)^{-1}\cdot e^{-s\boldsymbol{\sigma}}\\
&= \sqrt{\mathbf{Z}}\mathbf{Q}_1(\mathbf{I}+\mathbf{T}_1)(\mathbf{I}-\mathbf{T}_1)^{-1}\mathbf{Q}_1^{\mathrm{T}}\sqrt{\mathbf{Z}}^{-1}\\
&\quad - e^{-s\boldsymbol{\tau}}\cdot\sqrt{\mathbf{Z}}\mathbf{Q}_2(\mathbf{I}-\mathbf{T}_2)(\mathbf{I}+\mathbf{T}_2)^{-1}\mathbf{Q}_2^{\mathrm{T}}\sqrt{\mathbf{Z}}^{-1}\cdot e^{-s\boldsymbol{\sigma}}\\
&= \sqrt{\mathbf{Z}}\mathbf{K}\sqrt{\mathbf{Z}}^{-1}
\end{aligned}$$

Where,

(A.7) $\quad \mathbf{K} = \mathbf{Q}_1\boldsymbol{\Lambda}_1\mathbf{Q}_1^{\mathrm{T}} - e^{-s\boldsymbol{\tau}}\cdot\mathbf{Q}_2\boldsymbol{\Lambda}_2\mathbf{Q}_2^{\mathrm{T}}\cdot e^{-s\boldsymbol{\sigma}}$

$$\boldsymbol{\Lambda}_1 = (\mathbf{I}+\mathbf{T}_1)(\mathbf{I}-\mathbf{T}_1)^{-1}$$
$$\boldsymbol{\Lambda}_2 = (\mathbf{I}-\mathbf{T}_2)(\mathbf{I}+\mathbf{T}_2)^{-1}$$

Obviously, both $\boldsymbol{\Lambda}_1$ and $\boldsymbol{\Lambda}_2$ are diagonal matrices. All the diagonal elements of $\boldsymbol{\Lambda}_1$ are larger than 1, while all the diagonal elements of $\boldsymbol{\Lambda}_2$ are all less than 1. From above calculation, we have,

(A.8) $$\begin{aligned}\mathbf{I}-\mathbf{T} &= (\mathbf{I}+\mathbf{Z}\cdot\mathbf{A}_1)^{-1}\sqrt{\mathbf{Z}}\cdot\mathbf{K}\cdot\sqrt{\mathbf{Z}}^{-1}(\mathbf{I}-\mathbf{Z}\cdot\mathbf{A}_1)\\ (\mathbf{I}-\mathbf{T})^{-1} &= (\mathbf{I}-\mathbf{Z}\cdot\mathbf{A}_1)^{-1}\sqrt{\mathbf{Z}}\cdot\mathbf{K}^{-1}\cdot\sqrt{\mathbf{Z}}^{-1}(\mathbf{I}+\mathbf{Z}\cdot\mathbf{A}_1)\end{aligned}$$

Comparing (A.6) and (A.8), we can conclude, $s\mathbf{U}_1(s)$ has no pole in the region of the right half-plane plus the imaginary axis.

Now we prove, by reduction to absurdity, that $\mathbf{K}^{-1}$ has no pole in the region of the right half-plane plus imaginary axis. Assume $s_1$, one of the poles of $\mathbf{K}^{-1}$, is in the region of the right half-plane plus imaginary axis, then $\mathbf{K}(s_1)$ must be a singular complex matrix. Therefore, $\exists \boldsymbol{\varphi}\in\mathbb{C}^r$, $\|\boldsymbol{\varphi}\|=1$, and

$$\mathbf{K}(s_1)\boldsymbol{\varphi} = \mathbf{0}$$

Set $\mathbf{H}_1 = \mathbf{Q}_1\boldsymbol{\Lambda}_1\mathbf{Q}_1^{\mathrm{T}}$, $\mathbf{H}_2 = e^{-s_1\boldsymbol{\tau}}\cdot\mathbf{Q}_2\boldsymbol{\Lambda}_2\mathbf{Q}_2^{\mathrm{T}}\cdot e^{-s_1\boldsymbol{\sigma}}$. From (A.7) and the above equation, we have:

$$\mathbf{H}_1\boldsymbol{\varphi} = \mathbf{H}_2\boldsymbol{\varphi}$$

Consequently,

(A.9) $$\|\mathbf{H}_1\boldsymbol{\varphi}\| = \|\mathbf{H}_2\boldsymbol{\varphi}\|$$

$$\|\mathbf{H}_1\boldsymbol{\varphi}\| = \|\mathbf{Q}_1\boldsymbol{\Lambda}_1\mathbf{Q}_1^{\mathrm{T}}\boldsymbol{\varphi}\| = \sqrt{(\mathbf{Q}_1\boldsymbol{\Lambda}_1\mathbf{Q}_1^{\mathrm{T}}\boldsymbol{\varphi})^{\mathrm{H}}\mathbf{Q}_1\boldsymbol{\Lambda}_1\mathbf{Q}_1^{\mathrm{T}}\boldsymbol{\varphi}}$$
$$= \sqrt{\boldsymbol{\varphi}^{\mathrm{H}}\mathbf{Q}_1\boldsymbol{\Lambda}_1^2\mathbf{Q}_1^{\mathrm{T}}\boldsymbol{\varphi}} = \sqrt{(\mathbf{Q}_1^{\mathrm{T}}\boldsymbol{\varphi})^{\mathrm{H}}\boldsymbol{\Lambda}_1^2\mathbf{Q}_1^{\mathrm{T}}\boldsymbol{\varphi}} > \sqrt{r}$$

$$\|\mathbf{H}_2\boldsymbol{\varphi}\| = \|e^{-s_1\boldsymbol{\tau}}\cdot\mathbf{Q}_2\boldsymbol{\Lambda}_2\mathbf{Q}_2^{\mathrm{T}}\cdot e^{-s_1\boldsymbol{\sigma}}\boldsymbol{\varphi}\| \le \|e^{-s_1\boldsymbol{\tau}}\|_2\cdot\|e^{-s_1\boldsymbol{\sigma}}\|_2\cdot\|\mathbf{Q}_2\boldsymbol{\Lambda}_2\mathbf{Q}_2^{\mathrm{T}}\boldsymbol{\varphi}\|$$
$$\le \|e^{-s_1\boldsymbol{\tau}}\|_2\cdot\|e^{-s_1\boldsymbol{\sigma}}\|_2\cdot\sqrt{\boldsymbol{\varphi}^{\mathrm{H}}\mathbf{Q}_2\boldsymbol{\Lambda}_2^2\mathbf{Q}_2^{\mathrm{T}}\boldsymbol{\varphi}} < \|e^{-s_1\boldsymbol{\tau}}\|_2\cdot\|e^{-s_1\boldsymbol{\sigma}}\|_2\cdot\sqrt{r}$$

Because $\boldsymbol{\tau}$ and $\boldsymbol{\sigma}$ are positive vectors and $s_1$ is in the region of the right half-plane plus imaginary axis, we have,

$$\|e^{-s_1\boldsymbol{\tau}}\|_2 \le 1, \quad \|e^{-s_1\boldsymbol{\sigma}}\|_2 \le 1$$

As the result,

(A.10) $$\|\mathbf{H}_2\boldsymbol{\varphi}\| < \sqrt{r} < \|\mathbf{H}_1\boldsymbol{\varphi}\|$$

There is a contradiction between (A.9) and (A.10), thus we can conclude that $\mathbf{K}^{-1}$ has no pole in the region of the right half-plane plus imaginary axis. Further, from (A.8) and (A.6), we can conclude $s\mathbf{U}_1(s)$ and $(\mathbf{I}-\mathbf{T})^{-1}$ have no pole in the region of the right half-plane plus imaginary axis. Therefore we have,

$$\lim_{s \to 0}(\mathbf{I}-\mathbf{T})^{-1} s\mathbf{W}_1 = \lim_{s \to 0} \frac{s\mathbf{W}_1}{\mathbf{I}-\mathbf{T}}$$

$$= \frac{(\mathbf{I}+\mathbf{Z}\cdot\mathbf{A}_1)^{-1}\cdot(\mathbf{I}-\mathbf{Z}\cdot\mathbf{A}_2)\cdot(\mathbf{I}+\mathbf{Z}\cdot\mathbf{A}_2)^{-1}\cdot\mathbf{Z}\cdot\mathbf{b}_1 + (\mathbf{I}+\mathbf{Z}\cdot\mathbf{A}_1)^{-1}\mathbf{Z}\cdot\mathbf{b}_1}{\mathbf{I}-(\mathbf{I}+\mathbf{Z}\cdot\mathbf{A}_1)^{-1}\cdot(\mathbf{I}-\mathbf{Z}\cdot\mathbf{A}_2)\cdot(\mathbf{I}+\mathbf{Z}\cdot\mathbf{A}_2)^{-1}\cdot(\mathbf{I}-\mathbf{Z}\cdot\mathbf{A}_1)}$$

$$= \frac{(\mathbf{I}-\mathbf{Z}\cdot\mathbf{A}_2)\cdot(\mathbf{I}+\mathbf{Z}\cdot\mathbf{A}_2)^{-1}+\mathbf{I}}{(\mathbf{I}+\mathbf{Z}\cdot\mathbf{A}_1)-(\mathbf{I}-\mathbf{Z}\cdot\mathbf{A}_2)\cdot(\mathbf{I}+\mathbf{Z}\cdot\mathbf{A}_2)^{-1}\cdot(\mathbf{I}-\mathbf{Z}\cdot\mathbf{A}_1)} \mathbf{Z}\cdot\mathbf{b}_1$$

$$= \frac{\mathbf{I}+(\mathbf{I}-\mathbf{Z}\cdot\mathbf{A}_2)\cdot(\mathbf{I}+\mathbf{Z}\cdot\mathbf{A}_2)^{-1}}{\mathbf{I}-(\mathbf{I}-\mathbf{Z}\cdot\mathbf{A}_2)\cdot(\mathbf{I}+\mathbf{Z}\cdot\mathbf{A}_2)^{-1}+\left(\mathbf{I}+(\mathbf{I}-\mathbf{Z}\cdot\mathbf{A}_2)\cdot(\mathbf{I}+\mathbf{Z}\cdot\mathbf{A}_2)^{-1}\right)\mathbf{Z}\cdot\mathbf{A}_1} \mathbf{Z}\cdot\mathbf{b}_1$$

$$= \frac{\mathbf{I}}{\dfrac{\mathbf{I}-(\mathbf{I}-\mathbf{Z}\cdot\mathbf{A}_2)\cdot(\mathbf{I}+\mathbf{Z}\cdot\mathbf{A}_2)^{-1}}{\mathbf{I}+(\mathbf{I}-\mathbf{Z}\cdot\mathbf{A}_2)\cdot(\mathbf{I}+\mathbf{Z}\cdot\mathbf{A}_2)^{-1}}+\mathbf{Z}\cdot\mathbf{A}_1} \mathbf{Z}\cdot\mathbf{b}_1$$

$$= \frac{\mathbf{I}}{\dfrac{(\mathbf{I}+\mathbf{Z}\cdot\mathbf{A}_2)-(\mathbf{I}-\mathbf{Z}\cdot\mathbf{A}_2)}{(\mathbf{I}+\mathbf{Z}\cdot\mathbf{A}_2)+(\mathbf{I}-\mathbf{Z}\cdot\mathbf{A}_2)}+\mathbf{Z}\cdot\mathbf{A}_1} \mathbf{Z}\cdot\mathbf{b}_1$$

$$= \frac{\mathbf{I}}{\dfrac{2\mathbf{Z}\cdot\mathbf{A}_2}{2\mathbf{I}}+\mathbf{Z}\cdot\mathbf{A}_1} \mathbf{Z}\cdot\mathbf{b}_1$$

$$= \frac{\mathbf{I}}{\mathbf{Z}\cdot\mathbf{A}_1+\mathbf{Z}\cdot\mathbf{A}_2} \mathbf{Z}\cdot\mathbf{b}_1$$

$$= (\mathbf{A}_1+\mathbf{A}_2)^{-1}\cdot\mathbf{b}_1$$

Similarly, we get,

$$\lim_{s \to 0}(\mathbf{I}-\mathbf{T})^{-1} s\mathbf{W}_2 = \lim_{s \to 0} \frac{s\mathbf{W}_2}{\mathbf{I}-\mathbf{T}} = (\mathbf{A}_1+\mathbf{A}_2)^{-1}\mathbf{b}_2$$

According to the final-value theorem of Laplacian transformation, we have

$$\lim_{t \to \infty} u_1(t) = \lim_{s \to 0}(\mathbf{I}-\mathbf{T})^{-1}(s\mathbf{W}_1+s\mathbf{W}_2) = (\mathbf{A}_1+\mathbf{A}_2)^{-1}(\mathbf{b}_1+\mathbf{b}_2) = \mathbf{A}^{-1}\mathbf{b}$$

Similarly, we get,

$$\lim_{t \to \infty} u_2(t) = \mathbf{A}^{-1}\mathbf{b}$$

Therefore, Lemma A.1, the simple version of the convergence theorem, has been proved. The basic idea of this proof could be applied to prove the general convergence theorem of DTM.